\documentclass[12pt,a4paper]{article}

\usepackage{amsmath,amsfonts,amssymb,latexsym,graphics,epsfig,url}
\usepackage{xcolor}
\usepackage{amsthm}
\usepackage[english]{babel}
\usepackage{mathdots}
\usepackage{graphicx}
\usepackage{soul} 
\usepackage[utf8]{inputenc}
\usepackage{diagbox}
\usepackage[makeroom]{cancel}
\usepackage{multicol}
\usepackage{cancel}
\usepackage{float}              
\newtheorem{theorem}{Theorem}[section]

\newtheorem{lemma}[theorem]{Lemma}

\newtheorem{Example}[theorem]{Example}

\newtheorem{definition}[theorem]{Definition}

\DeclareMathOperator{\diam}{diam}
\DeclareMathOperator{\rad}{rad}
\DeclareMathOperator{\dist}{dist}

\DeclareMathOperator{\spec}{sp}

\def\x{\mbox{\boldmath $x$}}

\def\vec0{\mbox{\boldmath $0$}}

\def\L{\mbox{\boldmath $L$}}

\def\1{\mbox{\boldmath $1$}}

\begin{document}
	
\title{A note on an infinite family of graphs with all different integral Laplacian eigenvalues 
\thanks{This research has been supported by
AGAUR from the Catalan Government under project 2021SGR00434 and MICINN from the Spanish Government under project PID2020-115442RB-I00.
M. A. Fiol's research was also supported by a grant from the  Universitat Polit\`ecnica de Catalunya with references AGRUPS-2022 and AGRUPS-2023.}
}

\author{C. Dalf\'o$^a$ and M. A. Fiol$^b$\\
		{\small $^a$Dept. de Matem\`atica, Universitat de Lleida, 08700 Igualada (Barcelona), Catalonia}\\
		{\small {\tt cristina.dalfo@udl.cat}}\\
		{\small $^{b}$Dept. de Matem\`atiques, Universitat Polit\`ecnica de Catalunya, 08034  Barcelona, Catalonia} \\
		{\small Barcelona Graduate School of Mathematics} \\
		{\small  Institut de Matem\`atiques de la UPC-BarcelonaTech (IMTech)}\\
		{\small {\tt miguel.angel.fiol@upc.edu} }\\
	}

\date{}
\maketitle
 
\begin{abstract}
    In this note, we give an infinite family of optimal graphs called $G^+(d,c)$. They are optimal in the sense that they have the maximum possible number of vertices for given a diameter $d$ and the so-called `outer multiset dimension' $c$. We provide their spectra, which have the property that 
    their Laplacian eigenvalues are all different and integral. Finally, we also obtained their eigenvectors.
\end{abstract}



\section{Introduction }
\label{Intro}

Let $G$ be a finite, connected, undirected, and simple graph, with the vertex set $V(G)$ and the edge set $E(G)$. Let $W=\{w_1, \ldots, w_c\}$ be an ordered set of vertices in $G$ and let $u$ be a vertex of $G$. A {\em representation of $u$ with respect to $W$} is a vector $r(u|W)=(\dist(u,w_1), \ldots, \dist(u,w_c))$, where $\dist(u,w)=dist_G(u,w)$ is the distance between the vertices $u$ and $w$ in $G$. For every pair of distinct vertices $u,v$ in $V(G)$, if $r(u|W) \neq r(v|W)$, then $W$ is called {\em a resolving set for $G$}. For information on the multiset dimension of graphs, see Simanjuntak,  Vetrik, and Mulia \cite{SVM17}.

If we utilize a set of vertices $W=\{w_1,\ldots,w_c\}$ instead of an ordered set, then the representation of a vertex is a multiset instead of a vector. Let $u$ be a vertex of $G$. A {\em representation multiset of $u$ with respect to $W$}, denoted by $m(u|W)$, is a multiset of distances between $u$ and all vertices in $W$. For every pair of distinct vertices $u,v$ in $V(G)$, if $m(u|W)\neq m(v|W)$, then $W$ is called {\em a multiset resolving set for $G$}. To avoid an infinite multiset dimension, Gil-Pons, Ram\'irez-Cruz,  Trujillo-Rasua, and Yero \cite{GRTY19} (see also 
Klavzar, Kuziak, and Yero \cite{KKY23}) introduced the concept of an outer multiset dimension. Here, only vertices outside a multiset resolving set need to have distinct representation multiset. Formally, for every pair of distinct vertices $u,v$ in $V(G)\backslash W$, if $m(u|W)\neq m(v|W)$, then $W$ is called an {\em outer multiset resolving set for $G$}. 

In \cite{rados24}, Reyes, Araujo-Pardo, Dalf\'o, Olsen, and Simanjuntak gave the following new family of graphs, with some of its properties.

 \begin{definition}[\cite{rados24}]
\label{def:G(d,c)}
Given integers $d$ and $c$, the graph $G(d,c)$ has vertices identified with the combinations with repetition of $d$ elements $1,2,\ldots,d$ taken $c$ at a time.  Then, two vertices $x=x_1x_2\ldots x_c$ and $y=y_1y_2\ldots y_c$ are adjacent if and only if $|x_i-y_i|\in \{0,1\}$ for every $i=1,2,\ldots,c$.
\end{definition}
For instance, if $d=4$ and $c=3$, the vertex set of $G(4,3)$ contains the ${{6}\choose{3}}=20$ vertices $111$, $112$, $113$, 114, 122, 123, 124, 133, 134, 144, 222, 223, 224, 233, 234, 244, 333, 334, 344, and $444$ (notice that they can easily generated by lexicographic order). Thus, $G(4,3)$ is the graph induced by the black vertices in  Figure \ref{fig:G^+(3,3)} (right). Two more examples are  the graphs $G(2,6)$ and $G(3,3)$ induced by the black vertices in Figures \ref{fig:G^+(2,6)} (right) and \ref{fig:G^+(3,3)} (left), respectively.

\begin{figure}[t]
	\begin{center}
\includegraphics[width=12cm]{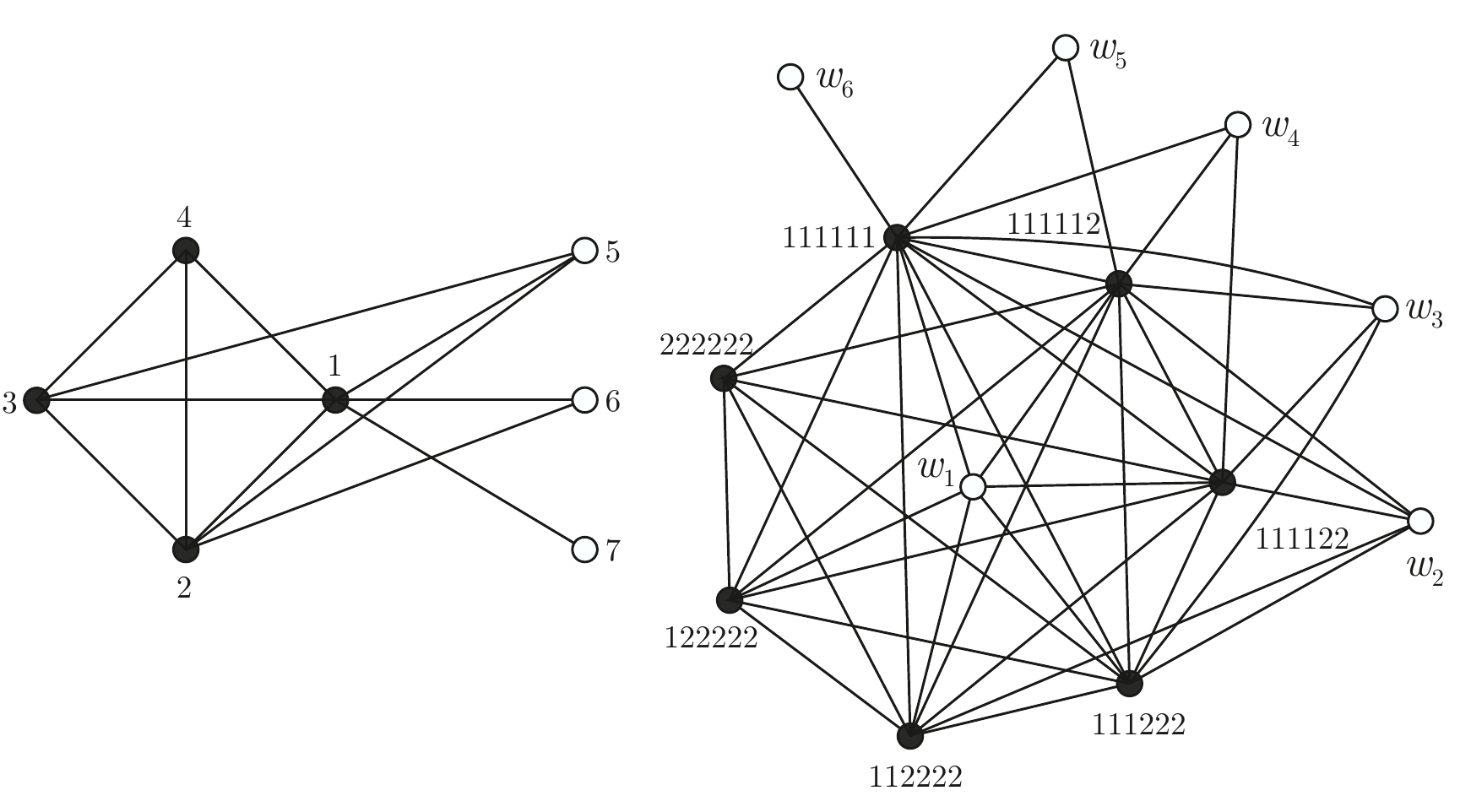}
  \vskip-.25cm
	\caption{The graphs $G^+(2,3)$ and $G^+(2,6)$, with diameter $d=2$, 
 and 
 7 and 13 
 vertices, respectively. The graphs induced by the black vertices are $G(2,3)$ and $G(2,6)$ with 
 4 and {7} vertices, respectively.}
		\label{fig:G^+(2,6)}
	\end{center}
\end{figure}

Some basic metric properties of the graph $G(d,c)$ are in the following result. 
\begin{lemma}[\cite{rados24}]
\label{lem:diam}
Let $x=x_1x_2\ldots x_c$ and $y=y_1y_2\ldots y_c$ two generic vertices of the graph $G(d,c)$. 
\begin{itemize}
\item[$(i)$]
   The distance between $x$ and $y$ is
   $$
   \dist(x,y)=\max_{i\in[1,c]}\{|x_i-y_i|\}.
   $$
\item[$(ii)$]
   The diameter of the graph $G(d,c)$ is 
   $$
   \diam(G(d,c))=d-1.
   $$
   \item[$(ii)$]
The radius of the graph $G(d,c)$ is 
   $$
   \rad (G(d,c))=\lfloor d/2\rfloor.
   $$
   \end{itemize}
\end{lemma}

From the family $G(d,c)$, they also defined another family of graphs $G^+(d,c)$.

\begin{definition}[\cite{rados24}]
Consider the graph $G^+(d,c)$ obtained from $G(d,c)$ by adding the vertex set $W=\{w_1,w_2,\ldots,w_c\}$, and connecting every vertex $w_i$ to all vertices of $G(d,c)$ with labels containing at least $i$ 1's for $i=1,2,\ldots,c$.
By Lemma \ref{lem:diam}, $G^+(d,c)$ has diameter $d$.
\end{definition}

The graphs $G^+(d,c)$ are optimal since they have the maximum possible number of vertices for given a diameter $d$ and an outer resolving multiset of $c$ vertices.

\begin{theorem}[\cite{rados24}]
\label{th:sharp}

Given positive integers $d$ and $c$, there exist a graph $G^+(d,c)$ with diameter $d$
containing an outer resolving multiset $W=\{w_1,w_2,\ldots,w_c\}$ with maximum number of vertices $n={d+c-1\choose d-1}+c$. 
The number of vertices of $G(d,c)$ is $n={d+c-1\choose d-1}$.
\end{theorem}



The graphs with diameter $2$ and $3$, $G^+(2,6)$ and  $G^+(3,3)$, are shown in Figures 
\ref{fig:G^+(2,6)} (right) and \ref{fig:G^+(3,3)} (left), respectively. 
Another example with diameter 4, $G^+(4,3)$, is shown in Figure  \ref{fig:G^+(3,3)} (left).

In the case of $d=2$, we can construct $G^+(2,c)$ (with its vertex labels) iteratively as follows:
\begin{enumerate}
\item
Start from $G^+(2,1)$, which is the path $x\sim y \sim w_1$ with labels  $x=2$ and $y=1$.
\item 
Given $G^+(2,c-1)$ with $n=2c-1$ vertices and resolving set $W_{c-1}=\{w_1,w_2,\ldots,$ $w_{c-1}\}$, relabel all is vertices adding a 2 at the end of their sequences, and add an isolated vertex with label  $w_c$.
\item 
Add another vertex with label $11\stackrel{(c)}{\ldots}1$ and join it to all the vertices of the previous graph $K_1\cup G^+(2,c-1)$.
\end{enumerate}
Thus, we proved that
\begin{equation}
    G^+(2,c)=K_1\cup(K_1\cup G^+(2,c-1)).
    \label{recurr}
\end{equation}

\begin{figure}[t]
	\begin{center}
		\includegraphics[width=14cm]{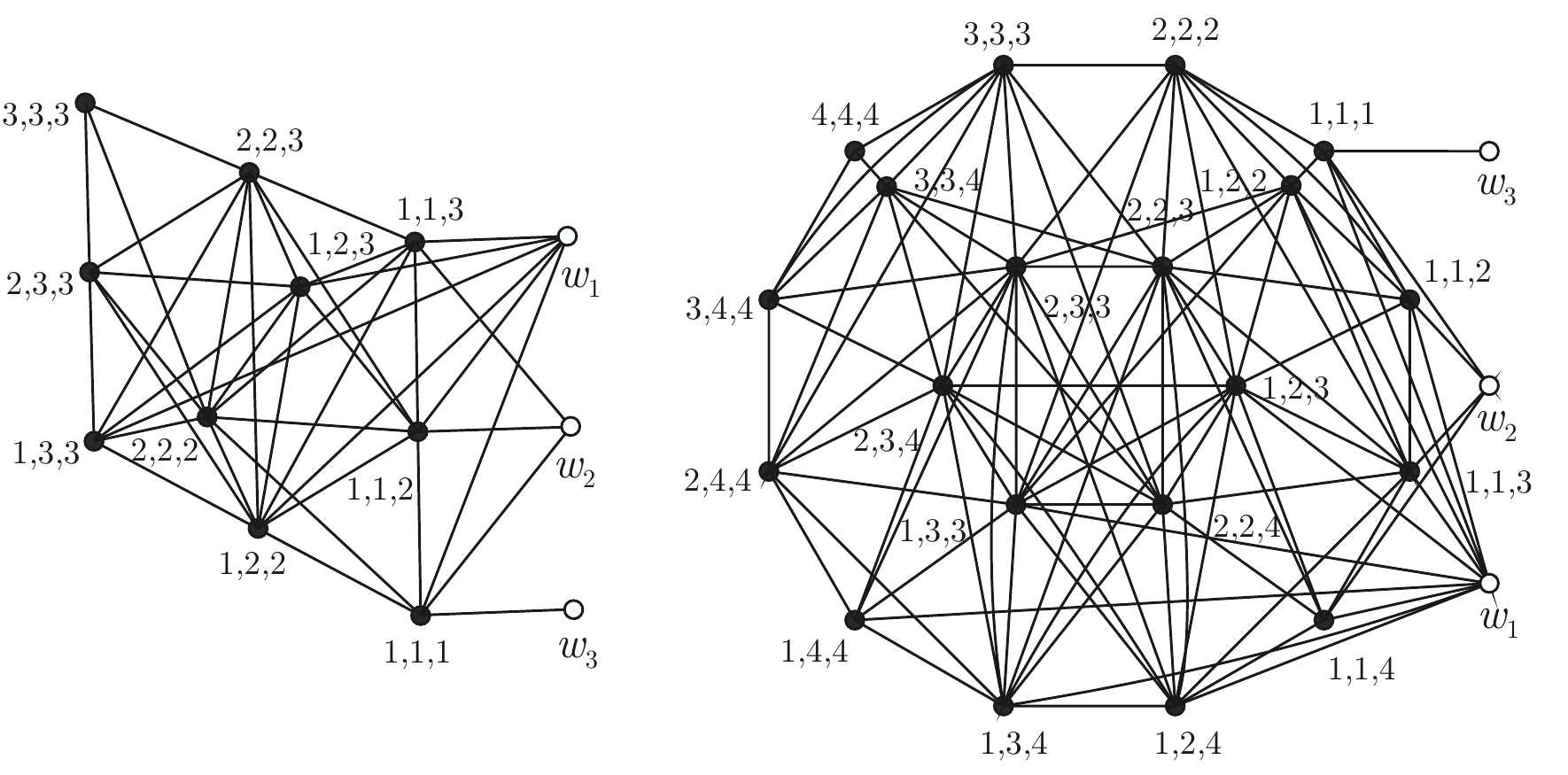}
  \vskip-.5cm
	\caption{The graphs $G^+(3,3)$ and $G^+(4,3)$ with 13 and 23 vertices, respectively. 
 The graphs induced by the black vertices are $G(3,3)$ and $G(4,3)$, with 
 10 and 20 vertices, respectively.}
		\label{fig:G^+(3,3)}
	\end{center}
\end{figure}


This note is structured as follows. In the current section, we gave the definitions and previous results on resolving sets, the outer multiset resolving set, and the two infinite families of graphs $G(d,c)$ and $G(d^+,c)$. In the following section, we obtain our main result, that is, the complete spectra and eigenspaces of $G^+(2,c)$.

\section{The Laplacian eigenvalues and eigenspaces of $G^+(2,c)$}

Fallat, Kirkland, Molitierno, and Neumann \cite{fkmn05} defined $S_{i,n}$ as 
the set of all integers from $0$ to $n$, excluding $i$. Then, they said that $S_{i,n}$ is {\em realizable} if there exists a graph $G$ on $n$ vertices with Laplacian matrix having $S_{i,n}$ as its spectrum. In this case, they also said that $G$ {\em realizes} $S_{i,n}$. Moreover, they characterized those $ i <n$ such that $S_{i,n}$ is Laplacian realizable and investigated the structure of the graphs that realize $S_{i,n}$. Concerning the latter, they prove the following (see \cite[Th. 2.3(c)]{fkmn05}).
\begin{theorem}[\cite{fkmn05}]
Let $G$ be a graph on $n\ge 6$ vertices. If $2\le i\le n-2$, then $G$ realizes $S_{i,n}$ if and only if $G=K_1 \vee (K_1 \cup H)$, where $H$ is a graph on $n-2$ vertices that realizes $S_{i-1,n-2}$.
\label{th:basic}
\end{theorem}


Using this notation, the following result shows that $G^+(2,c)$ realizes $S_{c+1,2c+1}$.

\begin{theorem}
Given any integer $c\ge 1$, the Laplacian spectrum of the graph $G^+(2,c)$ is
\begin{equation}
    \spec(G^+(2,c))=\{2c+1,2c,\ldots,c+2,c,c-1,\ldots,1,0\},
\end{equation}
with respective eigenvectors $\x^{(0)},\x^{(1)},\ldots,\x^{(2c)}$ being the columns of the following matrix:
\begin{equation}
\left(
\begin{array}{ccccc|c|ccccc|c}
-2c    & 0	         & \cdots  & 0	  & 0      & 0      & 0      & 0      & \cdots  & 0      & 0      & 1\\
1	   & -(2c-2) &     \cdots  & 0	  & 0      & 0      & 0      & 0      & \cdots  & 0      & -1	  & 1\\
1	   & 1	      & \cdots  & 0	  & 0      & 0      & 0      & 0      & \cdots  & -1	 & -1     & 1\\
\vdots & \vdots   & \ddots  & \vdots & \vdots & \vdots & \vdots & \vdots & \iddots & \vdots & \vdots & \vdots\\
1	   & 1		     & \cdots  & -4     & 0	   & 0      & 0      & -1     & \cdots  & -1	 & -1     & 1\\
1	   & 1	      & \cdots  & 1	  & -2     & 0      & -1     & -1     & \cdots  & -1	 & -1     & 1\\
1	   & 1	       & \cdots  & 1	  & 1	   & -1     & -1     & -1     & \cdots  & -1	 & -1     & 1\\
1	   & 1		      & \cdots  & 1	  & 1	   & 1      & -1     & -1     & \cdots  & -1	 & -1     & 1\\
1	   & 1		      & \cdots  & 1	  & 0	   & 0      & 3      & -1     & \cdots  & -1	 & -1     & 1\\
\vdots & \vdots    & \iddots & \vdots & \vdots & \vdots & \vdots & \vdots & \ddots  & \vdots & \vdots & \vdots\\
1	   & 1		      & \cdots  & 0	  & 0	   & 0      & 0      & 0      & \cdots  & -1	 & -1     & 1\\
1	   & 1		       & \cdots  & 0	  & 0	   & 0      & 0      & 0      & \cdots  & 2c-3   & -1     & 1\\
1	   & 0		      & \cdots  & 0	  & 0	   & 0      & 0      & 0      & \cdots  & 0      & 2c-1   & 1\\
\end{array}
\right).
\label{eigenvalues}
\end{equation}
\end{theorem}

\begin{proof}
Although we give below a direct proof by using eigenvectors, the spectrum of $G^+(2,c)$ follows from \eqref{recurr} and Theorem \ref{th:basic} of \cite{fkmn05}.
Indeed, first, it can be easily checked that, for $c=1,2,3$, $G^+(2,c)$ realizes $S_{c+1,2c+1}$. For instance, $G^+(2,3)$ in Figure \ref{fig:G^+(2,6)} has spectrum 
$\{7,6,5,3,2,1,0\}$. Then,  Theorem \ref{th:basic} can be applied recursively to prove the result. 

Concerning our direct proof involving the eigenvectors, we note first that,
by construction of the graph $G^+(2,c)$, its $(2c+1)\times (2c+1)$ Laplacian matrix is
$$
\L=\left(
\begin{array}{ccccccc|cccccc}
2c     & -1     & -1     & \cdots  & -1     & -1     & -1     & -1     & -1     & \cdots  & -1     & -1     & -1\\
-1     & 2c-1   & -1     & \cdots  & -1     & -1     & -1     & -1     & -1     & \cdots  & -1    & -1     & 0\\
-1     & -1     & 2c-2   & \cdots  & -1     & -1     & -1     & -1     & -1     & \cdots  & -1    & 0      & 0\\
\vdots & \vdots & \vdots & \ddots  & \vdots & \vdots & \vdots & \vdots & \vdots & \iddots & \vdots & \vdots & \vdots\\
-1     & -1     & -1     & \cdots  & c+2    & -1     & -1     & -1  & -1     & \cdots  & 0      & 0      & 0\\
-1     & -1     & -1     & \cdots  & -1    & c+1     & -1     & -1 & 0     & \cdots  & 0      & 0      & 0\\
-1     & -1     & -1     & \cdots  & -1     & -1     & c      & 0      & 0      & \cdots  & 0      & 0      & 0\\
\hline
-1     & -1     & -1     & \cdots  & -1     & -1     & 0      & c      & 0      & \cdots  & 0      & 0      & 0\\
-1     & -1     & -1     & \cdots  & -1     & 0      & 0      & 0      & c-1    & \cdots  & 0      & 0      & 0\\
\vdots & \vdots & \vdots & \iddots & \vdots & \vdots & \vdots & \vdots & \vdots & \ddots  & \vdots & \vdots & \vdots\\
-1     & -1     & -1      & \cdots  & 0     & 0      & 0      & 0      & 0      & \cdots  & 3      &  0     & 0\\
-1     & -1     & 0      & \cdots  & 0     & 0      & 0      & 0      & 0      & \cdots  & 0      &  2     & 0\\
-1     & 0      & 0      & \cdots  & 0     & 0      & 0      & 0      & 0      & \cdots  & 0      &  0     & 1\\
\end{array}
\right).
$$
The vertices corresponding to the entries of this Laplacian matrix are numbered first with the vertices of the graph $G(c,2)$, which is a complete graph on $c+1$ vertices and, finally, the vertices of the resolving set $w_1,\ldots,w_c$, giving a $2\times 2$ block matrix.


With this aim, first recall that, if $\x=(x_1,\ldots,x_n)^{\top}$ is an eigenvector of $\L$ with eigenvalue $\lambda$, then the Rayleigh quotient $\rho(\x)$ gives $\lambda$:
\begin{equation}
\label{Rayleigh}
\lambda=\rho(\x)=\frac{\x^{\top}\L\x}{\x^{\top}\x}=\frac{\sum_{ij\in E}|x_i-x_j|^2}{\sum_{i\in V} x_i^2}.
\end{equation}
Moreover, the set of edges of $G^+(d,c)=(V,E)$ can be partitioned as 
$E=\bigcup_{h=1}^{c}E_h$,
where 
$$
E_h=\{(h,j)\, : \, h+1\le j\le 2c+2-h\}\quad \mbox{for } h=1,2,\ldots,c.
$$

Hence, \eqref{Rayleigh} can be written as
$\rho(\x)=\frac{\sum_{h=1}^c N_h}{\sum_{i=1}^n x_i^2}$, where 
\begin{equation}
N_h=\sum_{hj\in E_h}|x_h-x_j|^2=\sum_{j=h+1}^{2c+2-h}|x_h-x_j|^2\quad \mbox{for } h=1,2,\ldots,c.
\label{Nh}
\end{equation}
With $n=2c+1$, Figure \ref{fig:squeme} shows the entries of $\x$ involved in the computation of each $N_h$ (notice the ``central symmetry'').

\begin{figure}[t]
	\begin{center}
		\includegraphics[width=10cm]{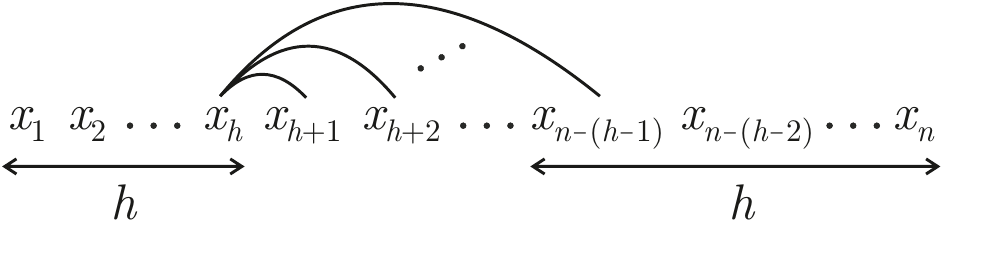}
  \vskip -.25cm
	\caption{Involved terms in the computation of $N_h$ (note that the curved lines correspond to the edges of $G^+(2,c)$ with endpoint $h$).}
		\label{fig:squeme}
	\end{center}
\end{figure}

Looking at the matrix \eqref{eigenvalues}, each eigenvector of $\L$ belongs to one of the following four classes (separated by the vertical lines in the matrix):
\begin{itemize}
\item[$(i)$] 
    For $r=0,\ldots,c-1$: 
   \begin{align*}
    \x^{(r)} &=(0,\stackrel{(r)}{\ldots},0,-2(c-r),1\stackrel{(2(c-r))}{\ldots\ldots}1, 0\stackrel{(r)}{\ldots}0)\\
     &=(0^r,-2(c-r), 1^{2(c-r)}, 0^r).
    \end{align*}
\item[$(ii)$] 
    For $r=c$: 
    $
    \x^{(c)}=(0^c,-1, 1,0^{c-1}).
    $
\item[$(iii)$] 
    For $r=c+1,\ldots,2c-1$: 
    $$
    \x^{(r)}=(0^{2c-r},-1^{2(r-c)+1}, 2(r-c)+1,0^{2c-r-1}).
    $$
\item[$(iv)$] 
    For $r=2c$: 
    $
    \x^{(2c)}=(1^{2c}).
    $
\end{itemize}
Now, we can compute the Rayleigh quotient by using the values of $N_h$ in \eqref{Nh} and the square norm $\|\x\|^2=\sum_{i\in V}x_i^2$:
\begin{itemize}
\item[$(i)$] 
    For $r=0,\ldots,c-1$: 
   \begin{align*}
    N_1 =\cdots =N_r &=[2(c-r)]^2+ 2(c-r)\quad \mbox{(provided that $r\neq 0$)};\\
    N_{r+1} &=[2(c-r)+1]^2\cdot2(c-r);\\
    N_{r+2}=\cdots =N_{c}&=0.
    \end{align*}
    Then,
    \begin{align*}
    \rho(\x^{(r)})&=\frac{rN_1+N_{r+1}}{\|\x^{(r)}\|^2}
    =\frac{r[[2(c-r)]^2+2(c-r)]+[2(c-r)+1]^2\cdot2(c-r)}{[2(c-r)]^2+2(c-r)}\\
    &= 2c-r+1 \in \{2c+1,2c,\ldots,c+2\}.
   \end{align*}
\item[$(ii)$] 
    For $r=c$:  $N_1 =\cdots =N_c =2$.
     Then,
    \begin{align*}
    \rho(\x^{(r)})&=\frac{cN_1}{\|\x^{(r)}\|^2}
    =\frac{2c}{2}= c.
   \end{align*}
\item[$(iii)$] 
    For $r=c+1,\ldots,2c-1$: 
    \begin{align*}
    N_1 =\cdots =N_{2c-r} &= 2(r-c)+1+[2(r-c)+1]^2;\\
    N_{2c-r+1}=\cdots =N_{c}&=0.
    \end{align*}
    Then,
    \begin{align*}
    \rho(\x^{(r)}) &=\frac{(2c-r)N_1}{\|\x^{(r)}\|^2}
   =\frac{(2c-r)[2(r-c)+1+[2(r-c)+1]^2]}{2(r-c)+1+[2(r-c)+1]^2}\\
    &= 2c-r \in \{c-1,c-2,\ldots,1\}.
   \end{align*}
\item[$(iv)$] 
    For $r=2c$: $N_1=\cdots=N_c=0$. Then, 
    $$
    \rho(\x^{(2c)})=\frac{cN_1}{\|\x^{(2c)}\|^2}=0.
    $$
\end{itemize}
Finally, since all the obtained eigenvalues are different, the corresponding eigenvectors must be linearly independent. This completes the proof.
\end{proof}

\begin{Example}
Let us show an example. The graph $G^+(2,3)=(V,E)$ on $n=7$ vertices in Figure \ref{fig:G^+(2,6)} (right) 
has Laplacian matrix
$$
\L=\left(
\begin{array}{cccc|ccc}
 6 & -1 & -1& -1& -1& -1& -1\\
-1 & 5 & -1& -1& -1& -1& 0  \\
-1 & -1& 4& -1& -1& 0& 0    \\
-1 & -1& -1& 3& 0& 0& 0     \\
\hline
-1 & -1& -1& 0& 3& 0& 0     \\
-1 & -1& 0& 0& 0& 2& 0      \\
-1 & 0 & 0& 0& 0& 0& 1
\end{array}
\right).
$$
We claim that the eigenvectors of $G^+(2,3)$ are the columns $\x^{(i)}$, for $i=0,1,\ldots,6$ of the following matrix:
$$
\left(
\begin{array}{ccc|c|cc|c}
-6 & 0	& 0  & 0  & 0      & 0  & 1\\
1 & -4 & 0  & 0  & 0        & {-}1  & 1\\
1 & 1 & -2  & 0  & {-}1  & {-}1  & 1\\
1 & 1 & 1  & -1  & {-}1  & {-}1  & 1\\
1 & 1 & 1  & 1  & {-}1  & {-}1  & 1\\
1 & 1 & 0  & 0  & {3}  & {-}1  & 1\\
1 & 0 & 0  & 0  & 0        & {5}  & 1
\end{array}
\right).
$$
Indeed, first notice that, with $\x=\x^{(i)}=(x_1,x_2,\ldots,x_7)^{\top}$ for some $i\in [0,6]$,
$$
\sum_{ij\in E}|x_i-x_j|^2=\sum_{h=1}^3 N_h=\sum_{j=2}^7|x_1-x_j|^2+\sum_{j=3}^6|x_2-x_j|^2+\sum_{j=4}^5|x_3-x_j|^2.
$$
\begin{itemize}
\item  
With $\x=\x^{(0)}=(-6,1,1,1,1,1,1)^{\top}$, we have $N_1=\sum_{j=2}^7|-6-1|^2$, and $N_2=N_3=0$. Then,
$$
\rho(\x)=\frac{\sum_{ij\in E}|x_1-x_j|^2}{\sum_{j\in V} x_i^2}=\frac{6\cdot 7^2}{6\cdot 7}=7.
$$
\item  
With $\x=\x^{(1)}=(0,-4,1,1,1,1,0)^{\top}$, we have $N_1=4^2+4\cdot 1^2=20$, $N_2=\sum_{j=3}^6 |-4-1|^2=100$, and $N_3=0$. Then,
$$
\rho(\x)=\frac{20+100}{\sum_{i\in V} x_i^2}
=\frac{120}{20}=6.
$$
\item  
With $\x=\x^{(2)}=(0,0,-2,1,1,0,0)^{\top}$, we have $N_1=(-2)^2+2\cdot 1^2=6$, $N_2=(-2)^2+\sum_{j=4}^6 1^2=6$, and $N_3=\sum_{j=4}^5 |-2-1|^2=18$. Then,
$$
\rho(\x)=\frac{6+6+18}{\sum_{i\in V} x_i^2}
=\frac{30}{6}=5.
$$
\item  
With $\x=\x^{({3})}=(0,0,0,-1,1,0,0)^{\top}$, we have $N_1=2$, $N_2=2$, and $N_3=2$. Then,
$$
\rho(\x)=\frac{2+2+2}{\sum_{i\in V} x_i^2}
=\frac{6}{2}=3.
$$
\item  
With $\x=\x^{({4})}=(0,0,-1,-1,-1,3,0)^{\top}$, we have $N_1=3+3^2=12$, $N_2=3+3^2=12$, and $N_3=0$. Then,
$$
\rho(\x)=\frac{12+12}{\sum_{i\in V} x_i^2}
=\frac{24}{12}=2.
$$

\item  
With $=\x=\x^{(5)}=(0,-1,-1,-1,{-}1,-1,5)^{\top}$, we have $N_1=5+5^2=30$, $N_2=0$, and $N_3=0$. Then,
$$
\rho(\x)=\frac{30}{\sum_{i\in V} x_i^2}
=\frac{30}{30}=1.
$$
\item  
With $\x=\x^{(6)}=(1,1,1,1,1,1,1)^{\top}$, we have $N_1=N_2=N_3=0$. Then,
$$
\rho(\x)=0.
$$
\end{itemize}
\end{Example}

\section*{Statements}

No data was used for the research described in the article. 

\noindent The authors have no conflict of interest.



\begin{thebibliography}{99}

\bibitem{fkmn05}
S. M. Fallat, S. J. Kirkland, J. J. Molitierno, and M. Neumann, 
On graphs whose Laplacian matrices have distinct integer eigenvalues, 
\textit{J. Graph Theory} \textbf{50} (2005) 162--174. 

\bibitem{GRTY19}
R. Gil-Pons, Y. Ramírez-Cruz, R. Trujillo-Rasua, and I. G. Yero,
Distance-based vertex identification in graphs: The outer multiset dimension,
\textit{App. Math. Comput.} \textbf{363} (2019) 124612.

\bibitem{KKY23}
S. Klavzar, D. Kuziak, and I. G. Yero,
Further contributions on the outer multiset dimension of graphs,
\textit{Results Math.} \textbf{78} (2023) 50.

\bibitem{rados24}
M. A. Reyes, G. M. Araujo-Pardo, C. Dalf\'o, M. Olsen, and R. Simanjuntak,
On the outer multiset dimensions of graphs,
submitted, 2024. 

\bibitem{SVM17}
R. Simanjuntak, T. Vetrik, and P. B. Mulia,
The multiset dimension of graphs, 
\texttt{arXiv:1711.00225}.



\end{thebibliography}
\end{document}